# APPROXIMATE WORD MATCHES BETWEEN TWO RANDOM SEQUENCES

By Conrad J. Burden, Miriam R. Kantorovitz[1] and Susan R. Wilson

*Australian National University, University of Illinois and Australian National University*

Given two sequences over a finite alphabet $\mathcal{L}$, the $D_2$ statistic is the number of $m$-letter word matches between the two sequences. This statistic is used in bioinformatics for expressed sequence tag database searches. Here we study a generalization of the $D_2$ statistic in the context of DNA sequences, under the assumption of strand symmetric Bernoulli text. For $k < m$, we look at the count of $m$-letter word matches with up to $k$ mismatches. For this statistic, we compute the expectation, give upper and lower bounds for the variance and prove its distribution is asymptotically normal.

**1. Introduction.** Methods for alignment-free sequence comparison are among the more recent tools being developed for sequence analysis in biology [16]. A disadvantage in the classical Smith–Waterman local alignment algorithm [13], as well as the popular search algorithms such as FASTA and BLAST, is that they assume conservation of contiguity between homologous segments. In particular, they overlook the occurrence of genetic shuffling [18]. Alignment-free sequence comparison methods are used to compensate for this problem.

A natural alignment-free comparison of two sequences is the number of $m$-letter word matches between the sequences. This statistic, called $D_2$, can be computed in linear time in the length of the sequences, which is also an advantage over the nonlinear local alignment algorithms. $D_2$ is used extensively for EST sequence database searches (e.g., [2, 3, 11] and in the software package STACK [6]).

In [10], Lippert, Huang and Waterman started a rigorous study of $D_2$ using the model of independent letters in DNA sequences. A formula for

Received May 2006; revised May 2007.
[1]Supported by Australian Research Council Discovery Grant DP0559260.
*AMS 2000 subject classifications.* 60F17, 92D20.
*Key words and phrases.* DNA sequences, sequence comparison, word matches.

the expectation was computed as well as upper and lower bounds for the variance. Limiting distributions, as the length of the sequences, $n$, and the size of the word, $m$, get large, were derived in some cases. The authors used Stein–Chen methods [5, 9, 14] to obtain the following results. When the underlying distribution of the alphabet is *nonuniform*, the distribution of $D_2$ has normal asymptotic behavior when $m/\log_b n < 1/6$. The logarithmic base $b$ is defined by $b = (\sum_{a\in\mathcal{L}} \xi_a^2)^{-1}$, where $\xi_a$ is the probability of a letter taking the value $a$. Following simulations, it was noted in [10] that the bound 1/6 above is too small. Our simulations in Section 6 suggest that the bound should be closer to 2.

Another asymptotic regime was identified in [10] when $m/\log_b n \geq 2$. In this case, the distribution of $D_2$ has compound Poisson asymptotic behavior. However, as pointed out in [17] and [1], the Poisson approximation is meaningful in this region only when $E(D_2)$ is not too small. To control this degenerate case, one needs to add the linear restriction $m = 2\log_b n + C$.

When the underlying distribution of the alphabet is *uniform*, it was proved in [9] that for $m = \alpha \log_b n + C$ with $0 < \alpha < 2$, the distribution of $D_2$ is also asymptotically normal.

A natural generalization of the $D_2$ statistic is to count the number of *approximate* $m$-word matches. For $k < m$, let $D_2^{(k)}$ be the number of $m$-word matches with up to $k$ mismatches between the two sequences. This statistic can be expressed in terms of a distance function. One can define the distance between two $m$-words to be the number of mismatches. A $k$-neighborhood of an $m$-word $\mathbf{w}$ is then all $m$-words that are at most $k$ distance from $\mathbf{w}$. The $D_2^{(k)}$ statistic is the number of $k$-neighborhood matches of $m$-words between two sequences.

In [12], Melko and Mushegian studied the $k$-distance and $k$-neighborhood match count between a probe of length $m$ and a random DNA sequence, under the assumption that the sequence is strand-symmetric Bernoulli text. They gave a formula for the expectation of the $k$-distance match count and the $k$-neighborhood match count. Melko and Mushegian suggested that methods of Lippert, Huang and Waterman in [10] could be used to obtain upper and lower bounds for the variance of $D_2^{(k)}$ and to analyze its asymptotic behavior.

In this paper we study the $D_2^{(k)}$ statistic under the strand-symmetric Bernoulli text assumption. We extend the method of [10] to give upper and lower bounds for the variance. We analyze the asymptotic behavior of the distribution of $D_2^{(k)}$ as $n$ and $m$ increase using the method of cumulants [8] rather than Stein's method. For $D_2$, the $k = 0$ case, this method improves the bound on $m/\log_b n$ obtained in [10] from 1/6 to 1/2.

The organization of this paper is as follows. In Section 2 we review definitions and introduce notation. In Section 3 we discuss the mean of $D_2^{(k)}$.

Section 4 is devoted to the variance of $D_2^{(k)}$. In Section 5 we prove normal asymptotic behavior of the distribution of $D_2^{(k)}$. Section 6 contains the results of numerical simulations, and a concluding summary is given in Section 7. A list of notations is provided at the end of Section 7.

**2. Preliminaries.** Let $\mathcal{L} = \{A, G, C, T\}$ with *strand-symmetric* probability measure $\xi = \{\xi_A, \xi_G, \xi_C, \xi_T\}$ and *perturbation parameter* $\eta$. That is, $-1 \leq \eta \leq 1$ is the unique number satisfying

$$\xi_A = \xi_T = \tfrac{1}{4}(1+\eta),$$
$$\xi_C = \xi_G = \tfrac{1}{4}(1-\eta).$$

Let $\mathbf{A} = A_1 A_2 \cdots A_n$ and $\mathbf{B} = B_1 B_2 \cdots B_n$ be two random sequences of length $n$ over the alphabet $\mathcal{L}$. We assume that $\mathbf{A}$ and $\mathbf{B}$ are Bernoulli texts, meaning, the letters (nucleotides) are independent and identically distributed (i.i.d.). We note that the assumption of both sequences having the same length is not essential for what follows and its main purpose is to simplify notation. Our results can be easily adapted to the case when the sequences are of different lengths.

Definition 2.1. Let $\mathbf{x}$ and $\mathbf{y}$ be two words of length $m$. We define the *distance* between $\mathbf{x}$ and $\mathbf{y}$ to be

$$\delta(\mathbf{x}, \mathbf{y}) = \text{ number of character mismatches between } \mathbf{x} \text{ and } \mathbf{y}.$$

For $k \leq m$, we say that $\mathbf{x}$ is a *$k$-distance match* of $\mathbf{y}$ if $\delta(\mathbf{x}, \mathbf{y}) = k$. When $\delta(\mathbf{x}, \mathbf{y}) \leq k$, then $\mathbf{x}$ is said to be a *$k$-neighbor* of $\mathbf{y}$.

Following the terminology and notation of [12], we have the following definition.

Definition 2.2. In the above setup, define the *perturbed binomial distribution* with perturbation parameter $\eta$ by

$$g_k(m, \eta, c) = h(m, \eta, c) u_k(m, \eta, c),$$

where $0 \leq c, k \leq m$ are integers and

$$h(m, \eta, c) = \frac{1}{4^m}(1-\eta)^c (1+\eta)^{m-c},$$

$$u_k(m, \eta, c) = \sum_{i=0}^{m-k} \binom{c}{i} \binom{m-c}{m-k-i} v_k(i, \eta, c),$$

$$v_k(i, \eta, c) = \left(\frac{3+\eta}{1-\eta}\right)^{c-i} \left(\frac{3-\eta}{1+\eta}\right)^{k-c+i}.$$

For an $m$-word $\mathbf{w}$ with $GC$-count $c_{\mathbf{w}}$, $h(m,\eta,c_{\mathbf{w}})$ is $\Pr(\mathbf{w})$, the probability of seeing $\mathbf{w}$. In the definition of $u_k(m,\eta,c)$, and in similar situations throughout the paper, we follow a general convention that $\binom{n}{a}=0$ if $a<0$ or $a>n$.

Note that when $\eta=0$, the perturbed binomial distribution is the binomial distribution with $g_k(m,0,c)=b_k(m,1/4)$, where

$$b_k(m,\rho)=\binom{m}{k}\rho^{m-k}(1-\rho)^k.$$

As observed in [12], if $\mathbf{T}$ is a strand-symmetric Bernoulli text of length $m$ and $\mathbf{q}$ is a (known) query text (=word) of length $m$, then the probability distribution of the distance $\delta(\mathbf{T},\mathbf{q})$ is a perturbed binomial distribution:

$$\Pr(\delta(\mathbf{T},\mathbf{q})=k)=g_k(m,\eta,c),$$

where $c$ is the $GC$-count in $\mathbf{q}$ and $\eta$ is the perturbation parameter of $\mathbf{T}$. Let

$$G_k(m,\eta,c)=\sum_{r=0}^{k} g_r(m,\eta,c)=\Pr(\delta(\mathbf{T},\mathbf{q})\le k)$$

be the cumulative distribution function of the distance.

2.1. *$k$-neighborhood matches.* Let $\mathbf{A}$ and $\mathbf{B}$ be two DNA sequences of length $n$. Assume the sequences are strand-symmetric Bernoulli text with perturbation parameter $\eta$. Let $0\le k\le m<n$ be integers.

Definition 2.3. Define the statistic $D_2^{(k)}=D(k,m,n)$ to be the number of $k$-neighborhood $m$-word matches between the sequences $\mathbf{A}$ and $\mathbf{B}$, including overlaps. Note that $D_2^{(0)}$ is the $D_2$ statistic of [10].

The $D_2^{(k)}$ statistic may be computed as follows.

Notation. For $1\le s\le t\le n$, write $\mathbf{A}[s,t]$ for the subsequence $A_sA_{s+1}\ldots A_t$.

Definition 2.4. Let $Y_{ij}^{(k)}$ be the $k$-neighborhood match indicator (starting) at position $(i,j)$ (position $i$ in sequence $\mathbf{A}$ and $j$ in $\mathbf{B}$). That is,

$$Y_{ij}^{(k)}=\begin{cases}1, & \text{if } \delta(\mathbf{A}[i,i+m-1],\mathbf{B}[j,j+m-1])\le k,\\ 0, & \text{otherwise.}\end{cases}$$

Then the $D_2^{(k)}$ statistic can be computed via

$$D_2^{(k)}=\sum_{(i,j)\in I} Y_{ij}^{(k)},$$

where the index set $I$ is

$$I = \{(i,j) \in \mathbb{N} \times \mathbb{N} : 1 \le i \le n-m+1, \text{ and } 1 \le j \le n-m+1\}.$$

For convenience, we write $\bar{n}$ for $n-m+1$.

## 3. The mean of $D_2^{(k)}$.

In this section we give a general formula for the mean of $D_2^{(k)}$ in terms of the perturbed binomial distribution and obtain estimates for it. The estimates will be used in later sections in order to prove normal asymptotic behavior of $D_2^{(k)}$.

First we compute the mean of $Y_{ij}^{(k)}$:

$$\begin{aligned} E[Y_{ij}^{(k)}] &= \Pr(Y_{ij}^{(k)} = 1) \\ &= \sum_{\mathbf{w} \in \mathcal{L}^m} \Pr(\delta(\mathbf{A}[i, i+m-1], \mathbf{w}) \le k) \Pr(\mathbf{B}[j, j+m-1] = \mathbf{w}) \\ &= \sum_{\mathbf{w} \in \mathcal{L}^m} G_k(m, \eta, c_{\mathbf{w}}) \Pr(\mathbf{w}), \end{aligned}$$

where $c_{\mathbf{w}}$ is the $GC$-count of $\mathbf{w}$.

From this we get formulas for the expectation of $D_2^{(k)}$:

$$\begin{aligned} E[D_2^{(k)}] &= \sum_{(i,j) \in I} E[Y_{ij}^{(k)}] \\ &= \bar{n}^2 \sum_{\mathbf{w} \in \mathcal{L}^m} \Pr(\mathbf{w}) G_k(m, \eta, c_{\mathbf{w}}). \end{aligned}$$

Remark 3.1. When $k = 0$, we have $G_k(m, \eta, c_{\mathbf{w}}) = \Pr(\mathbf{w})$ and

$$E[Y_{ij}^{(0)}] = \sum_{\mathbf{w} \in \mathcal{L}^m} (\Pr(\mathbf{w}))^2 = \left( \sum_{a \in \mathcal{L}} \xi_a^2 \right)^m.$$

This agrees with the formula given in Lippert, Huang and Waterman [10], $E[Y_{ij}] = p_2^m$, where $p_2 = \sum_{a \in \mathcal{L}} \xi_a^2$.

Definition 3.2. For $t > 1$, let

$$p_t = \sum_{a \in \mathcal{L}} \xi_a^t.$$

Remark 3.3. Note that $p_t = E[(\xi_X)^{t-1}]$. Hence, by the Cauchy–Schwarz inequality,

$$p_3 \ge p_2^2,$$

where equality holds if and only if the distribution is uniform: $\xi_a = 1/|\mathcal{L}|$.

3.1. *Estimates.* The purpose of these estimates is to explain the asymptotic behavior of $D_2^{(k)}$, rather than to provide a computational tool. Hence, these estimates are by no means optimal.

First we estimate the function $g_k(m,\eta,c_{\mathbf{w}})$. Without loss of generality, assume $\eta>0$. From Definition 2.2, we have upper bounds:

$$v_k(i,\eta,c) \le \left(\frac{3+\eta}{1-\eta}\right)^k,$$

$$u_k(m,\eta,c) \le \left(\frac{3+\eta}{1-\eta}\right)^k \sum_{i=0}^{m-k} \binom{c}{i}\binom{m-c}{m-k-i} = \left(\frac{3+\eta}{1-\eta}\right)^k \binom{m}{k}.$$

Remembering that $h(m,\eta,c_{\mathbf{w}}) = \Pr(\mathbf{w})$ and using similar estimates for the lower bound we get

$$\Pr(\mathbf{w})\binom{m}{k}\left(\frac{3-\eta}{1+\eta}\right)^k \le g_k(m,\eta,c_{\mathbf{w}}) \le \Pr(\mathbf{w})\binom{m}{k}\left(\frac{3+\eta}{1-\eta}\right)^k$$

and

$$\Pr(\mathbf{w})\sum_{r=0}^{k}\binom{m}{r}\left(\frac{3-\eta}{1+\eta}\right)^r \le G_k(m,\eta,c_{\mathbf{w}}) \le \Pr(\mathbf{w})\sum_{r=0}^{k}\binom{m}{r}\left(\frac{3+\eta}{1-\eta}\right)^r. \tag{1}$$

Hence, for $E[Y_{ij}^{(k)}] = \sum_{\mathbf{w}\in\mathcal{L}^m}\Pr(\mathbf{w})G_k(m,\eta,c_{\mathbf{w}})$, we have

$$p_2^m\sum_{r=0}^{k}\binom{m}{r}\left(\frac{3-\eta}{1+\eta}\right)^r \le E[Y_{ij}^{(k)}] \le p_2^m\sum_{r=0}^{k}\binom{m}{r}\left(\frac{3+\eta}{1-\eta}\right)^r. \tag{2}$$

Finally, since $E[D_2^{(k)}] = \sum_{i,j=1}^{\bar{n}} E[Y_{ij}^{(k)}]$, we have

$$\bar{n}^2 p_2^m\sum_{r=0}^{k}\binom{m}{r}\left(\frac{3-\eta}{1+\eta}\right)^r \le E[D_2^{(k)}] \le \bar{n}^2 p_2^m\sum_{r=0}^{k}\binom{m}{r}\left(\frac{3+\eta}{1-\eta}\right)^r. \tag{3}$$

Remark 3.4. For $k=0$ (the exact matches case), (3) gives $E[D_2^{(0)}] = \bar{n}^2 p_2^m$, which agrees with the expectation computed in [10].

Remark 3.5. When $\eta=0$ (the uniform case), $p_2=\frac{1}{4}$ and the upper and lower bounds in (3) are equal. Hence,

$$E[D_2^{(k)}] = \bar{n}^2\frac{1}{4^m}\sum_{r=0}^{k}\binom{m}{r}3^r$$

$$= \bar{n}^2 \sum_{r=0}^{k} \binom{m}{r} \left(\frac{3}{4}\right)^r \left(\frac{1}{4}\right)^{m-r}$$

$$= \bar{n}^2 \sum_{r=0}^{k} b_r(m, 1/4).$$

**4. The variance.** In this section we give a lower and upper bound for the variance of $D_2^{(k)}$. The lower bound is used later to prove asymptotic normality of $D_2^{(k)}$. The upper bound is not optimal, but is comparable with that given in [10] for the $k = 0$ case. We start this section by stating the main results in Propositions 4.1 and 4.2. We then prove several technical lemmas and finish with the proofs of Propositions 4.1 and 4.2.

PROPOSITION 4.1. $\mathrm{Var}(D_2^{(k)}) \geq \ell = \ell(n, m, k)$, *where*

$$\ell = \bar{n}^2 \left[ (2\bar{n} - 4m + 2)\,(m-1)\,k^2 \left(\frac{3-\eta}{1+\eta}\right)^{2k} p_2^{2m} \left(\frac{p_3}{p_2^2} - 1\right) \right] + O(n^2 m^{k+2} p_2^m).$$

PROPOSITION 4.2.

$$\begin{aligned} \mathrm{Var}(D_2^{(k)}) &\leq \bar{n}^2 (2\bar{n} - 4m + 2) m^{2k} \left(\frac{3+\eta}{1-\eta}\right)^{2k} \left[ 2p_3 \left(\frac{1-p_3^m}{1-p_3}\right) - p_3^m \right] \\ &\quad - \bar{n}^2 (2\bar{n} - 4m + 2) \left[ 2p_2^2 \left(\frac{1-p_2^{2m}}{1-p_2^2}\right) - p_2^{2m} \right] \\ &\quad + \bar{n}^2 (2m-1)^2 p_2^m \sum_{r=0}^{k} \binom{m}{r} \left(\frac{3+\eta}{1-\eta}\right)^r . \end{aligned} \tag{4}$$

To compute the variance of $D_2^{(k)} = \sum_{(i,j)\in I} Y_{ij}^{(k)}$, we need to compute the covariances $\mathrm{Cov}(Y_{ij}^{(k)}, Y_{i'j'}^{(k)})$. For this, we use techniques from [17]. To shorten the indices' notation, let $u = (i, j)$ and $v = (i', j')$.

In the following definition we use notation and terminology from [17], Chapter 11.

DEFINITION 4.3. Let $J_u = \{v = (i', j') : |i' - i| < m \text{ or } |j' - j| < m\}$. Then $J_u$ is the *dependency neighborhood* of $Y_u^{(k)}$ in the sense that if $v \notin J_u$, then $Y_u^{(k)}$ and $Y_v^{(k)}$ are independent. The dependency neighborhood can be decomposed into two parts, *accordion* and *crabgrass*, $J_u = J_u^a \cup J_u^c$, where

$$J_u^a = \{v = (i', j') \in J_u : |i' - i| < m \text{ and } |j' - j| < m\} \quad \text{and} \quad J_u^c = J_u \setminus J_u^a.$$

Let $u \in I$. When $v \notin J_u$, $\mathrm{Cov}(Y_u^{(k)}, Y_v^{(k)}) = 0$. To estimate $\mathrm{Cov}(Y_u^{(k)}, Y_v^{(k)})$ when $v \in J_u$, we look at the two cases: $v \in J_u^c$ (crabgrass) and $v \in J_u^a$ (accordion). We will see that crabgrasses contribute the dominant term of $\mathrm{Var}(D_2^{(k)})$ in the cases we are interested in, that is, for $m = O(\log n)$. Hence for accordions we only give a crude approximation of $\mathrm{Cov}(Y_u^{(k)}, Y_v^{(k)})$. We start by proving the following positivity lemma.

LEMMA 4.4. (i) *For $v \in J_u^c$, $Y_u^{(k)}$ and $Y_v^{(k)}$ are nonnegatively correlated. That is,*

$$\mathrm{Cov}(Y_u^{(k)}, Y_v^{(k)}) \geq 0.$$

(ii) *For $v$ in the main diagonal of $J_u^a$,* $\mathrm{Cov}(Y_u^{(k)}, Y_v^{(k)}) \geq 0$.

PROOF. We will use the following notation. For $r \geq 0$ define

$$\begin{aligned}
Y_{ij}^{(k)}(r) &= \text{the } k\text{-neighbor match indicator between two } r\text{-words at } (i,j). \\
\Delta_1(r) &= \delta(A[i, i+r-1], B[j, j+r-1]) \\
&= \text{number of mismatches in an } r\text{-word match at } (i,j). \\
\Delta_2(r) &= \delta(A[i'+m-r, i'+m-1], B[j'+m-r, j'+m-1]) \\
&= \text{number of mismatches in an } r\text{-word match} \\
&\quad\ \text{at } (i'+m-r, j'+m-r).
\end{aligned}$$

To prove part (i), let $u = (i,j)$, $v = (i',j') \in J_u^c$. Write $t = i' - i$ and $s = j' - j$. By symmetry, we may assume $v$ is in the first quadrant of $J_u^c$, that is, $0 \leq t \leq m-1$, and $m \leq s$. We have

$$\begin{aligned}
&E[Y_u^{(k)} Y_v^{(k)}] \\
&\quad = \Pr(Y_u^{(k)} = 1, Y_v^{(k)} = 1) \\
&\quad = \sum_{l_1, l_2 = 0}^{t} \Pr(\Delta_1(t) = l_1) \Pr(\Delta_2(t) = l_2) \\
&\qquad\qquad \times \Pr(Y_{(i', j+t)}^{(k-l_1)}(m-t) = 1, Y_{(i',j')}^{(k-l_2)}(m-t) = 1) \\
&\quad = \sum_{l_1, l_2 = 0}^{t} \Pr(\Delta_1(t) = l_1) \Pr(\Delta_2(t) = l_2) \qquad (5) \\
&\qquad\qquad \times \left[ \sum_{\mathbf{w} \in \mathcal{L}^{m-t}} \Pr(\mathbf{w}) G_{k-l_1}(m-t, \eta, c_{\mathbf{w}}) G_{k-l_2}(m-t, \eta, c_{\mathbf{w}}) \right]
\end{aligned}$$

$$= \sum_{\mathbf{w}\in\mathcal{L}^{m-t}} \Pr(\mathbf{w}) \left[\sum_{l_1} \Pr(\Delta_1(t)=l_1) G_{k-l_1}(m-t,\eta,c_{\mathbf{w}})\right] \times \left[\sum_{l_2} \Pr(\Delta_2(t)=l_2) G_{k-l_2}(m-t,\eta,c_{\mathbf{w}})\right].$$

For $\mathbf{w}\in\mathcal{L}^{m-t}$, let

$$f_t(\mathbf{w}) = \sum_{l} \Pr(\Delta(t)=l) G_{k-l}(m-t,\eta,c_{\mathbf{w}}), \tag{6}$$

where $\Delta(t)$ is the distance between two random $t$-words. Then (5) says that

$$E[Y_u^{(k)} Y_v^{(k)}] = E[(f_t(W))^2].$$

Similarly we get

$$E[Y_u^{(k)}]E[Y_v^{(k)}] = E[f_t(W)]^2.$$

Hence,

$$\mathrm{Cov}(Y_u^{(k)}, Y_v^{(k)}) = \mathrm{Var}(f_t(W)) \geq 0. \tag{7}$$

For part (ii), let $u=(i,j)$ and $v=(i',j')\in J_u^a$'s main diagonal, that is, $v=(i+t,j+t)$ with $|t|\leq m-1$. By symmetry, we may assume $0\leq t$. As before, let $\Delta_1(t)$ be the number of mismatches in a $t$-word match at $(i,j)$, and let $\Delta_2(t)$ be the number of mismatches at $(i+m,j+m)$. Then

$$\begin{aligned} E[Y_u^{(k)} Y_v^{(k)}] &= \sum_{l_1,l_2=0}^{t} \Pr(\Delta_1(t)=l_1)\Pr(\Delta_2(t)=l_2) \\ &\qquad \times \Pr(Y_v^{(k-l_1)}(m-t)=1, Y_v^{(k-l_2)}(m-t)=1) \\ &= \sum_{l_1,l_2=0}^{t} \Pr(\Delta_1(t)=l_1)\Pr(\Delta_2(t)=l_2) \\ &\qquad \times \Pr(Y_v^{(\min\{k-l_1,k-l_2\})}(m-t)=1) \end{aligned}$$

and

$$\begin{aligned} E[Y_u^{(k)}]E[Y_v^{(k)}] &= \sum_{l_1,l_2=0}^{t} \Pr(\Delta_1(t)=l_1)\Pr(\Delta_2(t)=l_2) \\ &\qquad \times \Pr(Y_v^{(k-l_1)}(m-t)=1)\Pr(Y_v^{(k-l_2)}(m-t)=1). \end{aligned}$$

Since

$$\begin{aligned} &\Pr(Y_v^{(\min\{k-l_1,k-l_2\})}(m-t)=1) \\ &\qquad \geq \Pr(Y_v^{(k-l_1)}(m-t)=1)\Pr(Y_v^{(k-l_2)}(m-t)=1), \end{aligned}$$

we have that $\mathrm{Cov}(Y_u^{(k)}, Y_v^{(k)}) \geq 0$. $\square$

Remark 4.5. From (7) we have that for $v \in J_u^c$, $\mathrm{Cov}(Y_u^{(k)}, Y_v^{(k)}) = \mathrm{Var}(f_t(W))$ for appropriate $t$. When computing $f_t(\mathbf{w})$ in (6), it is worth noting the following:

1. $\Pr(\Delta(t) = l) = \sum_{\mathbf{x}\in\mathcal{L}^t} \Pr(\mathbf{x}) g_l(t, \eta, c_{\mathbf{x}}) = \sum_{\mathbf{x}\in\mathcal{L}^t} \Pr(\mathbf{x})^2 u_l(t, \eta, c_{\mathbf{x}})$. From Section 3.1 we also have
$$\binom{t}{l}\left(\frac{3-\eta}{1+\eta}\right)^l \le u_l(t, \eta, c_{\mathbf{x}}) \le \binom{t}{l}\left(\frac{3+\eta}{1-\eta}\right)^l.$$
Hence
$$p_2^t \binom{t}{l}\left(\frac{3-\eta}{1+\eta}\right)^l \le \Pr(\Delta(t) = l) \le p_2^t \binom{t}{l}\left(\frac{3+\eta}{1-\eta}\right)^l.$$
2. For $k - l \ge m - t$ and $\mathbf{w} \in \mathcal{L}^{m-t}$, $G_{k-l}(m - t, \eta, c_{\mathbf{w}}) = 1$.

Remark 4.6. When $k = 0$ (exact matches case), and $v \in J_u^c$ with $t$ as above, we get:

For $\mathbf{w} \in \mathcal{L}^{m-t}$,
$$f_t(\mathbf{w}) = \Pr(\Delta(t) = 0) G_0(m - t, \eta, c_{\mathbf{w}}) = p_2^t \Pr(\mathbf{w}).$$
Since $E[\Pr(W)] = \sum_{\mathbf{w}\in\mathcal{L}^{m-t}} \Pr(\mathbf{w})^2 = p_2^{m-t}$ and $E[\Pr(W)^2] = \sum_{\mathbf{w}\in\mathcal{L}^{m-t}} \Pr(\mathbf{w})^3 = p_3^{m-t}$, we have that
$$\mathrm{Cov}(Y_u^{(k)}, Y_v^{(k)}) = \mathrm{Var}(f_t(W)) = p_2^{2t}\,\mathrm{Var}(\Pr(\mathbf{w})) = p_2^{2t}(p_3^{m-t} - p_2^{2(m-t)}).$$
This agrees with the computations in [17], Section 11.5.2.

Remark 4.7. When $\eta = 0$ (uniform case), $f_t(\mathbf{w})$ does not depend on $\mathbf{w}$ and hence $\mathrm{Var}(f_t(W)) = 0$. Therefore, in this case, for $v \in J_u^c$, $\mathrm{Cov}(Y_u^{(k)}, Y_v^{(k)}) = 0$.

Next we look at the following special crabgrass case.

Lemma 4.8. *Let $u = (i, j)$ and $v = (i', j') \in J_u^c$ with $t = |i' - i| = m - 1$ or (by symmetry) $|j' - j| = m - 1$. Then*
$$\mathrm{Cov}(Y_u^{(k)}, Y_v^{(k)}) = [\Pr(\Delta(m-1) = k)]^2 (p_3 - p_2^2),$$
*where $\Delta(m - 1)$ is the distance between two random $(m - 1)$-words.*

*Hence, by Remark 4.5,*
$$\left[\binom{m-1}{k}\left(\frac{3-\eta}{1+\eta}\right)^k\right]^2 p_2^{2m}\left(\frac{p_3}{p_2^2} - 1\right) \le \mathrm{Cov}(Y_u^{(k)}, Y_v^{(k)}) \le \left[\binom{m-1}{k}\left(\frac{3+\eta}{1-\eta}\right)^k\right]^2 p_2^{2m}\left(\frac{p_3}{p_2^2} - 1\right).$$

PROOF. By (7),

$$\mathrm{Cov}(Y_u^{(k)}, Y_v^{(k)}) = \mathrm{Var}(f_t(W)).$$

Let $\mathbf{w} \in \mathcal{L}$. Then, using Remark 4.5,

$$\begin{aligned} f_t(\mathbf{w}) &= \sum_{l=0}^{k} \Pr(\Delta(m-1) = l) G_{k-l}(1, \eta, c_{\mathbf{w}}) \\ &= \Pr(\Delta(m-1) = k) G_0(1, \eta, c_{\mathbf{w}}) + \sum_{l=0}^{k-1} \Pr(\Delta(m-1) = l) \cdot 1 \\ &= \Pr(\Delta(m-1) = k) \xi_{\mathbf{w}} + \sum_{l=0}^{k-1} \Pr(\Delta(m-1) = l). \end{aligned}$$

Note that $\Pr(\Delta(m-1) = k)$ and $\sum_{l=0}^{k-1} \Pr(\Delta(m-1) = l)$ do not depend on $\mathbf{w}$. Hence,

$$\begin{aligned} \mathrm{Var}(f_t(W)) &= \mathrm{Var}\left( \Pr(\Delta(m-1) = k) \xi_W + \sum_{l=0}^{k-1} \Pr(\Delta(m-1) = l) \right) \\ &= [\Pr(\Delta(m-1) = k)]^2 \mathrm{Var}(\xi_W). \end{aligned}$$

As noted before (Remark 4.6, with $m - t = 1$), $\mathrm{Var}(\xi_W) = p_3 - p_2^2$. □

For the accordion case, we use the following crude estimate.

LEMMA 4.9. *For* $\mathbf{u}, \mathbf{v} \in I$, $\mathrm{Cov}(Y_u^{(k)}, Y_v^{(k)}) = O(p_2^m m^k)$.

PROOF.

$$\begin{aligned} |\mathrm{Cov}(Y_u^{(k)}, Y_v^{(k)})| &\le \sqrt{\mathrm{Var}(Y_u^{(k)}) \mathrm{Var}(Y_v^{(k)})} \\ &= \mathrm{Var}(Y_u^{(k)}) \le E[(Y_u^{(k)})^2] = E[Y_u^{(k)}] \\ &\le p_2^m \sum_{r=0}^{k} \binom{m}{r} \left( \frac{3+\eta}{1-\eta} \right)^r \qquad \text{from (2)} \\ &= O(p_2^m m^k). \end{aligned} \tag{8}$$

□

4.1. *Proof of Proposition* 4.1. We now prove the lower bound formula for $\mathrm{Var}(D_2^{(k)})$ stated in Proposition 4.1.

PROOF OF PROPOSITION 4.1. First we split the variance into the contributions of crabgrasses and accordions:

$$\begin{aligned}\mathrm{Var}(D_2^{(k)}) &= \sum_{u,v\in I} \mathrm{Cov}(Y_u^{(k)}, Y_v^{(k)}) \\ &= \sum_{u\in I}\sum_{v\in J_u^c} \mathrm{Cov}(Y_u^{(k)}, Y_v^{(k)}) + \sum_{u\in I}\sum_{v\in J_u^a} \mathrm{Cov}(Y_u^{(k)}, Y_v^{(k)}).\end{aligned}$$

Next we look at the crabgrasses:

$$\begin{aligned}&\sum_{u\in I}\sum_{v\in J_u^c} \mathrm{Cov}(Y_u^{(k)}, Y_v^{(k)}) \\ &\quad\geq \sum_{u=(i,j)\in I} \sum_{\substack{v=(i',j')\in J_u^c \\ |i'-i|=m-1 \text{ or } |j'-j|=m-1}} \mathrm{Cov}(Y_u^{(k)}, Y_v^{(k)}) \qquad \text{by Lemma 4.4} \\ &\quad\geq \sum_{u=(i,j)\in I} \sum_{\substack{v=(i',j')\in J_u^c \\ |i'-i|=m-1 \text{ or } |j'-j|=m-1}} \binom{m-1}{k}^2 \left(\frac{3-\eta}{1+\eta}\right)^{2k} p_2^{2m}\left(\frac{p_3}{p_2^2}-1\right) \\ &\qquad\qquad \text{by Lemma 4.8} \\ &\quad= 2\bar{n}^2\left[(2\bar{n}-4m+2)\binom{m-1}{k}^2 \left(\frac{3-\eta}{1+\eta}\right)^{2k} p_2^{2m}\left(\frac{p_3}{p_2^2}-1\right)\right].\end{aligned}$$

Finally we consider the contribution of the accordions to the variance:

$$\begin{aligned}\left|\sum_{u\in I}\sum_{v\in J_u^a} \mathrm{Cov}(Y_u^{(k)}, Y_v^{(k)})\right| &\leq \sum_{u\in I}\sum_{v\in J_u^a} |\mathrm{Cov}(Y_u^{(k)}, Y_v^{(k)})| \\ &= \sum_{u\in I}\sum_{v\in J_u^a} O(p_2^m m^k) \qquad \text{by Lemma 4.9} \\ &= \bar{n}^2(2m-1)^2 O(p_2^m m^k) = O(n^2 m^{k+2} p_2^m).\end{aligned}$$

Then

$$\begin{aligned}\mathrm{Var}(D_2^{(k)}) &= \sum_{u\in I}\sum_{v\in J_u^c} \mathrm{Cov}(Y_u^{(k)}, Y_v^{(k)}) + O(n^2 m^{k+2} p_2^m) \\ &\geq 2\bar{n}^2\left[(2\bar{n}-4m+2)\binom{m-1}{k}^2 \left(\frac{3-\eta}{1+\eta}\right)^{2k} p_2^{2m}\left(\frac{p_3}{p_2^2}-1\right)\right] \\ &\quad + O(n^2 m^{k+2} p_2^m). \qquad \square\end{aligned}$$

4.2. *Proof of Proposition 4.2.* The first two terms in (4) come from crabgrasses. Let $u=(i,j)\in I$, $v=(i',j')\in J_u^c$ with $i'=i+t$, $0\leq t\leq m-1$. We need to bound $\mathrm{Cov}(Y_u^{(k)}, Y_v^{(k)}) = \mathrm{Var}(f_t(W))$:

$$\begin{aligned}
E[f_t(W)^2] &= \sum_{\mathbf{w}\in\mathcal{L}^{m-t}} \Pr(\mathbf{w}) \left[\sum_{l=0}^{m-1} \Pr(\Delta(t)=l) G_{k-l}(m-t,\eta,c_{\mathbf{w}})\right]^2 \\
&\le \sum_{\mathbf{w}\in\mathcal{L}^{m-t}} \Pr(\mathbf{w})^3 \left[\sum_{l=0}^{m-1} \Pr(\Delta(t)=l) \sum_{r=0}^{k-l} \binom{m-t}{r} \left(\frac{3+\eta}{1-\eta}\right)^r\right]^2 \\
&\qquad\qquad \text{by (1)} \\
&\le \sum_{\mathbf{w}\in\mathcal{L}^{m-t}} \Pr(\mathbf{w})^3 \left[(m-t)^k \left(\frac{3+\eta}{1-\eta}\right)^k\right]^2 \left[\sum_{l=0}^{m-1} \Pr(\Delta(t)=l)\right]^2 \\
&= \left[(m-t)^k \left(\frac{3+\eta}{1-\eta}\right)^k\right]^2 \sum_{\mathbf{w}\in\mathcal{L}^{m-t}} \Pr(\mathbf{w})^3 \\
&= \left[(m-t)^k \left(\frac{3+\eta}{1-\eta}\right)^k\right]^2 p_3^{m-t} \\
&\le p_3^{m-t} m^{2k} \left(\frac{3+\eta}{1-\eta}\right)^{2k},
\end{aligned}$$

and similarly,

$$\begin{aligned}
E[f_t(W)]^2 &= \left[\sum_{\mathbf{w}\in\mathcal{L}^{m-t}} \Pr(\mathbf{w}) \sum_{l=0}^{m-1} \Pr(\Delta(t)=l) G_{k-l}(m-t,\eta,c_{\mathbf{w}})\right]^2 \\
&\ge \left[\sum_{\mathbf{w}\in\mathcal{L}^{m-t}} \Pr(\mathbf{w})^2 \sum_{l=0}^{m-1} \Pr(\Delta(t)=l) \sum_{r=0}^{k-l} \binom{m-t}{r} \left(\frac{3-\eta}{1+\eta}\right)^r\right]^2 \\
&\ge \left[\sum_{\mathbf{w}\in\mathcal{L}^{m-t}} \Pr(\mathbf{w})^2 \sum_{l=0}^{m-1} \Pr(\Delta(t)=l)\right]^2 \\
&= \left[\sum_{\mathbf{w}\in\mathcal{L}^{m-t}} \Pr(\mathbf{w})^2\right]^2 = [p_2^{m-t}]^2.
\end{aligned}$$

Hence,

$$\mathrm{Cov}(Y_u^{(k)}, Y_v^{(k)}) \le p_3^{m-t} m^{2k} \left(\frac{3+\eta}{1-\eta}\right)^{2k} - p_2^{2(m-t)}.$$

Summing up over all $u$'s and $v$'s and using

$$\sum_{u\in I} \sum_{v\in J_u^c} q^{m-t} = \bar{n}^2(2\bar{n} - 4m + 2)\left[2q\left(\frac{1-q^m}{1-q}\right) - q^m\right],$$

with $q = p_3$ and $p_2^2$, respectively, yields

$$\sum_{u\in I}\sum_{v\in J_u^c} \operatorname{Cov}(Y_u^{(k)}, Y_v^{(k)}) \le \bar{n}^2(2\bar{n} - 4m + 2)m^{2k}\left(\frac{3+\eta}{1-\eta}\right)^{2k}\left[2p_3\left(\frac{1-p_3^m}{1-p_3}\right) - p_3^m\right]$$
$$- \bar{n}^2(2\bar{n} - 4m + 2)\left[2p_2^2\left(\frac{1-p_2^{2m}}{1-p_2^2}\right) - p_2^{2m}\right].$$

The last term in (4) comes from accordions:

$$\begin{aligned}\sum_{u\in I}\sum_{v\in J_u^a} \operatorname{Cov}(Y_u^{(k)}, Y_v^{(k)}) &\le \sum_{u\in I}\sum_{v\in J_u^a} |\operatorname{Cov}(Y_u^{(k)}, Y_v^{(k)})| \\ &\le \sum_{u\in I}\sum_{v\in J_u^a} p_2^m \sum_{r=0}^{k} \binom{m}{r}\left(\frac{3+\eta}{1-\eta}\right)^r \qquad \text{from (8)} \\ &= \bar{n}^2(2m-1)^2 p_2^m \sum_{r=0}^{k}\binom{m}{r}\left(\frac{3+\eta}{1-\eta}\right)^r.\end{aligned}$$

**5. Asymptotic behavior.** We will need the following central limit theorem of Janson [8] for certain sums of dependent random variables. To state it, we first recall the definition of dependency graph.

DEFINITION 5.1. A graph $\Gamma$ is a dependency graph for a family of random variables if the following holds:

1. There is a one-to-one correspondence between the random variables and the vertices of the graph.
2. If $V_1$ and $V_2$ are two disjoint sets of vertices of $\Gamma$ such that there is no edge $(v_1, v_2)$ in $\Gamma$ with $v_1 \in V_1$ and $v_2 \in V_2$, then the corresponding sets of random variables are independent.

Also recall that the maximal degree of a graph is the maximal number of edges attached to a single vertex.

THEOREM 5.2 ([8], Theorem 2). *Suppose that for each $n$, $\{W_{ni}\}_{i=1}^{N_n}$ is a family of bounded random variables; $|W_{ni}| \le C_n$ almost surely. Suppose further that $\Gamma_n$ is a dependency graph for this family and let $M_n$ be the maximal degree of $\Gamma_n$ (if $\Gamma_n$ has no edges, set $M_n = 1$). Let $S_n = \sum_{i=1}^{N_n} W_{ni}$ and $\sigma_n^2 = \operatorname{Var}(S_n)$. If there exists an integer $t$ such that*

$$(N_n/M_n)^{1/t} M_n C_n/\sigma_n \to 0 \qquad \text{as } n \to \infty, \tag{9}$$

*then*

$$(S_n - E(S_n))/\sigma_n \xrightarrow{d} \mathcal{N}(0,1) \qquad \text{as } n \to \infty.$$

Next we state and prove our main theorem.

THEOREM 5.3. *Assume that the four letters of the alphabet $\mathcal{L}$ are not uniformly distributed. That is, the perturbation parameter $\eta$ is not zero. Let $\mu_n = E[D_2^{(k)}]$ and $\sigma_n = \sqrt{\mathrm{Var}(D_2^{(k)})}$.*

*For $m = \alpha \log_{1/p_2}(n) + C$ with $0 \le \alpha < 1/2$ and $C$ a constant, and fixed $k$ such that $0 \le k < m$,*

$$\frac{D_2^{(k)} - \mu_n}{\sigma_n} \overset{d}{\Longrightarrow} \mathcal{N}(0,1) \qquad \text{as } n \to \infty.$$

PROOF. We apply Theorem 5.2 to the match indicator random variables $Y_{ij}^{(k)}$. In this case, the dependency graph has $\bar{n}^2$ vertices and edges may be defined by connecting the vertex $(i,j)$ with $(i',j')$ if $|i'-i| < m$ or $|j'-j| < m$. Hence, in the notation of Theorem 5.2, $N_n = \bar{n}^2$; $C_n = 1$; the maximal degree of $\Gamma_n$ is the size of a dependency neighborhood:

$$M_n = |J_u| = (2m-1)(2\bar{n} - 2m + 1);$$

and $S_n = D_2^{(k)}$.

Let $m = \alpha \log_{1/p_2}(n) + C$ with $0 \le \alpha$ (and $k$ fixed). Then for $\alpha < 1$, the lower bound, $\ell$, for $\sigma_n^2$ in Proposition 4.1, has the property

$$\begin{aligned}\ell &\ge \bar{n}^2 \left[(2\bar{n} - 4m + 2)\left(\frac{3-\eta}{1+\eta}\right)^{2k} p_2^{2m}\left(\frac{p_3}{p_2^2} - 1\right)\right] + O(n^2 m^{k+2} p_2^m) \\ &= C_1 n^{3-2\alpha} + O(n^{2-\alpha}(\log(n))^{k+2}) \qquad \text{where } C_1 > 0 \text{ is a constant} \\ &\sim n^{3-2\alpha} \qquad \text{since } \alpha < 1.\end{aligned}$$

Therefore, in condition (9) we have

$$\begin{aligned}&\frac{(N_n/M_n)^{1/t} M_n C_n}{\sigma_n} \\ &\quad = \frac{(\bar{n}^2/((2m-1)(2\bar{n}-2m+1)))^{1/t}(2m-1)(2\bar{n}-2m+1)}{\sigma_n} \\ &\quad \le \left(\frac{\bar{n}^2}{(2m-1)(2\bar{n}-2m+1)}\right)^{1/t} (2m-1)(2\bar{n}-2m+1) \\ &\qquad \times \left(\left\{\bar{n}^2\left[(2\bar{n}-4m+2)\left(\frac{3-\eta}{1+\eta}\right)^{2k} p_2^{2m}\left(\frac{p_3}{p_2^2}-1\right)\right]\right.\right. \\ &\qquad\qquad \left.\left. + O(n^2 m^{k+2} p_2^m)\right\}^{1/2}\right)^{-1}\end{aligned} \tag{10}$$

$$\sim \frac{(\alpha \log_{1/p_2}(n))^{1-1/t} n^{1+1/t}}{[n^{3-2\alpha}]^{1/2}} = \frac{(\alpha \log_{1/p_2}(n))^{1-1/t}}{n^{1/2-1/t-\alpha}}$$

$$\to 0 \qquad \text{as } n \to \infty, \quad \text{if } 1/2 - 1/t - \alpha > 0.$$

Thus, for $\alpha < 1/2$, we can find $t$ large enough such that $1/2 - 1/t - \alpha > 0$. □

In [10], Lippert, Huang and Waterman used a variation on Stein's result ([15], page 110), due to Dembo and Rinott ([7], Theorem 4.2), to prove the following result, under the assumption of i.i.d. letters, for the $D_2 = D_2^{(0)}$ statistic. Let $\mathcal{L}$ be an alphabet set of size $|\mathcal{L}| > 1$ with nonuniform probability measure $\xi$. Then for $m = \alpha \log_{1/p_2}(n) + C$ with $0 \leq \alpha < 1/6$,

$$\frac{D_2(n) - \mu_n}{\sigma_n} \overset{d}{\Longrightarrow} \mathcal{N}(0,1) \qquad \text{as } n \to \infty.$$

Following simulations, it was noted in [10] that the bound 1/6 above is too small. Our simulations in Section 6 suggest that the bound should be closer to 2.

By adjusting the proof of Theorem 5.3 to an alphabet set of any size $|\mathcal{L}| > 1$, we can improve the bound on $\alpha$ from 1/6 to 1/2. Thus we have the following theorem.

THEOREM 5.4. *Let $\mathcal{L}$ be an alphabet set of size $|\mathcal{L}| > 1$ with nonuniform probability measure $\xi$. Then for $m = \alpha \log_{1/p_2}(n) + C$ with $0 \leq \alpha < 1/2$,*

$$\frac{D_2(n) - \mu_n}{\sigma_n} \overset{d}{\Longrightarrow} \mathcal{N}(0,1) \qquad \textit{as } n \to \infty.$$

PROOF. From the lower bound for $D_2$ in [10],

$$\begin{aligned} \operatorname{Var}(D_2) \geq \bar{n}^2 \bigg[ &(2m-1)(2\bar{n} - 4m + 2) p_2^{2m} (p_3/p_2^2 - 1) \\ &+ p_2^m \left( \frac{1 + p_2 - 2p_2^2}{1 - p_2} - (2m-1) p_2^m \right) \bigg] \\ \sim \alpha \log_{1/p_2}(n)(n^{3-2\alpha}) \qquad &\text{since } (p_3/p_2^2 - 1) > 0 \text{ by Remark 3.3.} \end{aligned}$$

Hence, in (10) in the proof of Theorem 5.3 we now have

$$\begin{aligned} &\frac{(N_n/M_n)^{1/t} M_n C_n}{\sigma_n} \\ &\quad = \frac{(\bar{n}^2/((2m-1)(2\bar{n} - 2m + 1)))^{1/t} (2m-1)(2\bar{n} - 2m + 1)}{\sigma_n} \end{aligned}$$

$$
\begin{aligned}
&\le \left(\frac{\bar{n}^2}{(2m-1)(2\bar{n}-2m+1)}\right)^{1/t}(2m-1)(2\bar{n}-2m+1)\\
&\quad\times\left(\bar{n}\left[(2m-1)(2\bar{n}-4m+2)p_2^{2m}\left(\frac{p_3}{p_2^2}-1\right)\right.\right.\\
&\qquad\qquad\left.\left.+\,p_2^m\left(\frac{1+p_2-2p_2^2}{1-p_2}-(2m-1)p_2^m\right)\right]^{1/2}\right)^{-1}\\
&\sim \frac{(\alpha\log_{1/p_2}(n))^{1-1/t}n^{1+1/t}}{[\alpha\log_{1/p_2}(n)n^{3-2\alpha}]^{1/2}} = \frac{(\alpha\log_{1/p_2}(n))^{1/2-1/t}}{n^{1/2-1/t-\alpha}}\\
&\to 0 \qquad \text{as } n\to\infty,\ \text{if } 1/2-1/t-\alpha>0.
\end{aligned}
$$

The rest of the proof is the same as the proof of Theorem 5.3. □

REMARK 5.5. When the underlying distribution is uniform, $\mathrm{Var}(D_2^{(k)}) \sim n^2p_2^{2m}$. Hence, this method of proof fails to show normal asymptotic behavior in the uniform case. In fact, for $k=0$, the distribution of $D_2^{(0)}$ is *not* normal when $|\mathcal{L}|=2, m=1$ and $n\to\infty$ (see [10]).

**6. Numerical simulations.** We have carried out numerical simulations of pairs of randomly generated sequences of length $n = 100\times 2^i$, $i=0,\dots,4$ with the nonuniform letter distribution $\xi_A=\xi_T=\frac{1}{3}$, $\xi_C=\xi_G=\frac{1}{6}$. The statistic $D_2^{(k)}$ was calculated for each sequence pair using an algorithm based on that given in [12]. Kolmogorov–Smirnov $p$-values [4] for the standardized statistic $(D_2^{(k)}-\mu_n)/\sigma_n$ compared with the standard normal distribution for sample sizes of 2500 sequence pairs are shown in Figure 1. Samples of $D_2^{(k)}$ which are a close approximation to the normal distribution will have $p$-values distributed uniformly on the interval $[0,1]$, whereas samples which are a poor approximation to the normal distribution will have small $p$-values.

Entries in the tables in Figure 1 are shaded to indicate proximity of samples to the standard normal distribution, with lighter shades signaling a better agreement. The white diagonal line in each table is $m = 2\log_{1/p_2} n +$ const. The numerical evidence strongly suggests that

$$\frac{D_2^{(k)}-\mu_n}{\sigma_n} \xrightarrow{d} \mathcal{N}(0,1) \qquad \text{as } n\to\infty,$$

where the limit is taken along any line $m=\alpha\log_{1/p_2}(n)+C$ with $0\le\alpha<2$ for fixed $k$ and $C$.

**7. Conclusions.** We have studied the $D_2^{(k)}$ statistic as suggested in [12], and defined it as the number of $m$-word matches with up to $k$ mismatches between two sequences of length $n$, for strand-symmetric Bernoulli texts with a nonuniform letter distribution. We have extended methods applied in [10] to give upper and lower bounds on the variance, and have also studied the asymptotic behavior as the sequence length $n$ and word length $m$ tend to infinity for fixed $k$.

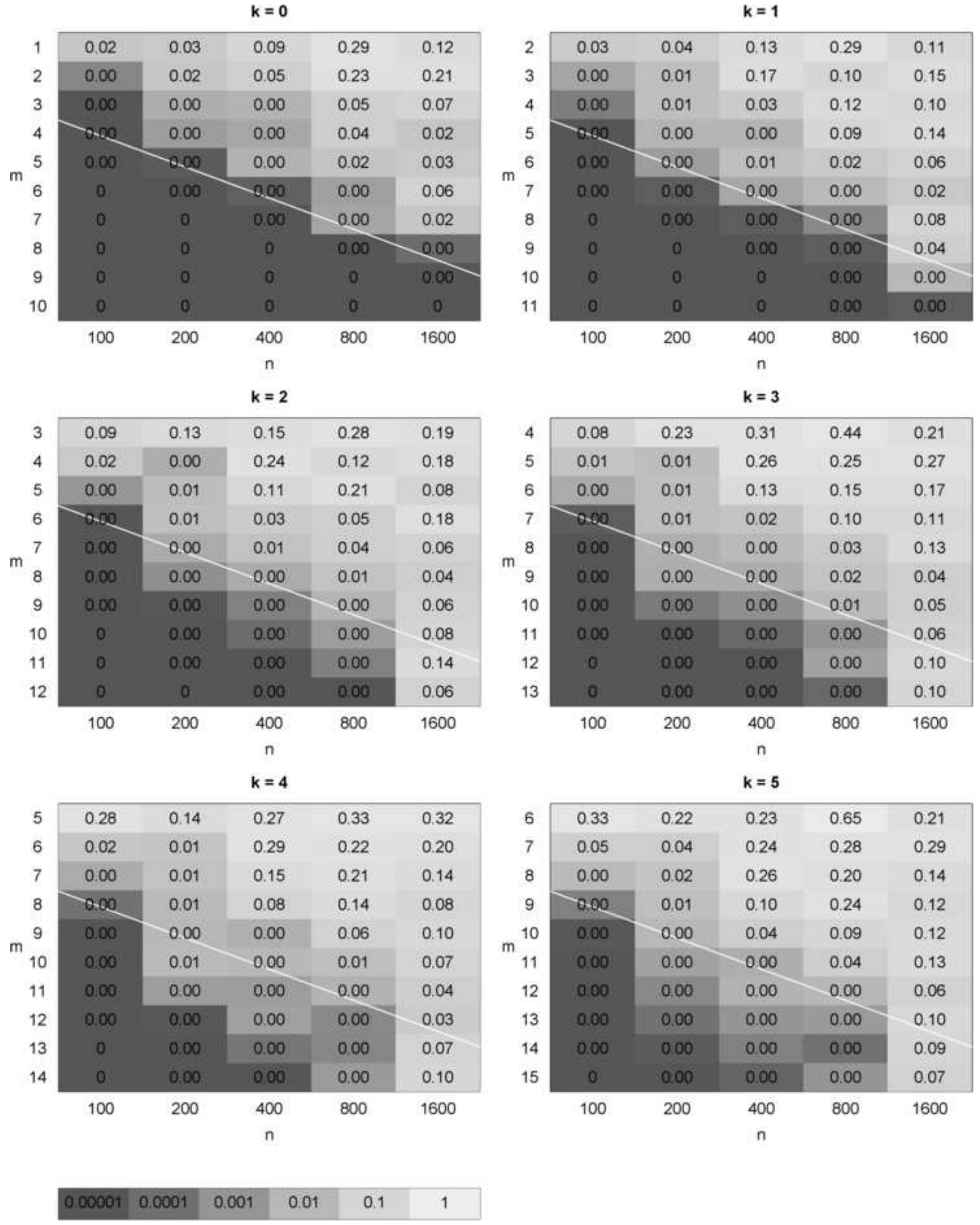

FIG. 1. *Kolmogorov–Smirnov p-values for nonuniform $D_2^{(k)}$ with letter distribution $\xi_A = \xi_T = \frac{1}{3}$, $\xi_C = \xi_G = \frac{1}{6}$ compared with a normal distribution. The white diagonal lines are $m = \alpha \log_{1/p_2} n + const.$, with $\alpha = 2$ and $1/p_2 = 1/\sum_{a\in\{A,C,G,T\}} \xi_a^2 = \frac{18}{5}$.*

We have proved that the asymptotic distributional behavior of $D_2^{(k)}$ is normal as $n \to \infty$ for a pair of strand-symmetric Bernoulli texts provided the limit is taken along any line $m = \alpha \log_{1/p_2} n + C$ with $0 \leq \alpha < \frac{1}{2}$. For $k = 0$ this result is also shown to hold for a pair of Bernoulli texts with any nonuniform letter distribution. This improves the previous bound for the $k = 0$ case of $\alpha < \frac{1}{6}$ given in [10].

We have also carried out numerical simulations of strand-symmetric texts with letter distribution $\xi_A = \xi_T = \frac{1}{3}$, $\xi_C = \xi_G = \frac{1}{6}$. These simulations strongly suggest that the optimum restriction on asymptotic normal behavior may be as high as $\alpha < 2$. This is consistent with simulations in [10] and their result that the asymptotic distributional behavior of $D_2^{(0)}$ is a compound Poisson distribution for $\alpha \geq 2$.

**List of notations.**

$D_2$: The number of $m$-letter word matches between two given sequences.

$D_2^{(k)}$: The number of $m$-letter word matches with up to $k$ mismatches ($0 \leq k \leq m$) between two given sequences.

$g_k(m, \eta, c)$: The perturbed binomial distribution (Definition 2.2). Given a strand-symmetric Bernoulli text of length $m$ and perturbation parameter $\eta$, and an $m$-word with $GC$-content $c$, this is the probability distribution of the number of character mismatches between the text and the $m$-word.

$G_k(m, \eta, c)$: The cumulative distribution function of the perturbed binomial distribution $g_k(m, \eta, c)$.

$h_k(m, \eta, c)$: For a given $m$-word with $GC$-content $c$, the probability that the word occurs at a given site in a strand-symmetric Bernoulli string with perturbation parameter $\eta$. (See Definition 2.2.)

$J_u$: The dependency neighborhood of $Y_u^{(k)}$, where $u = (i, j)$; that is, the word locations $v = (i', j')$ such that either the word at $i'$ overlaps the word at $i$ in the first sequence, or the word at $j'$ overlaps the word at $j$ in the first sequence, or both. (See Definition 4.3.)

$J_u^a$: The accordion, that is, the subset of $J_u$ such that both the word at $i'$ overlaps the word at $i$ in the first sequence, and the word at $j'$ overlaps the word at $j$ in the first sequence.

$J_u^c$: The crabgrass, that is, the subset of $J_u$ such that either the word at $i'$ overlaps the word at $i$ in the first sequence, or the word at $j'$ overlaps the word at $j$ in the first sequence, but not both.

$k$: The number of mismatches.
$m$: The word length.
$n$: The length of each of the two sequences.
$\bar{n}$: $n-m+1$, the number of possible locations of an $m$-word in a sequence of length $n$.
$p_t$: $\sum_{a\in\mathcal{L}} \xi_a^t$, where the sum is taken over the alphabet $\mathcal{L}$.
$u_k(m,\eta,c), v_k(i,\eta,c)$: Functions occurring in the definition of $g_k(m,\eta,c)$. (See Definition 2.2.)
$Y_u^{(k)}$ or $Y_{ij}^{(k)}$: The indicator random variable for the event that the $m$-word starting at position $i$ in the first sequence has no more than $k$ mismatches with the $m$-word starting at position $j$ in the second sequence. We use the convention $u=(i,j), v=(i',j')$ throughout.
$\eta$: The perturbation parameter for a strand-symmetric Bernoulli text. (See Section 2.)
$\xi_a$: The probability of finding the letter $a$ at a given location in a strand-symmetric Bernoulli string.

C. J. Burden
John Curtin School of Medical Research
and Mathematical Sciences Institute
Australian National University
Canberra, A.C.T. 0200
Australia
E-mail: conrad.burden@anu.edu.au

M. R. Kantorovitz
Department of Mathematics
University of Illinois
Urbana, Illinois 61801
USA
E-mail: ruth@math.uiuc.edu

S. R. Wilson
Mathematical Sciences Institute
Australian National University
Canberra, A.C.T. 0200
Australia
E-mail: sue.wilson@anu.edu.au